\documentclass[12pt,english,a4paper]{article}

\usepackage[english]{babel}

\usepackage{amssymb}

\usepackage[T1]{fontenc}

\usepackage{enumerate}

\usepackage{amscd}

\usepackage{amsthm}

\usepackage{amsmath}

\usepackage{array}

\usepackage{delarray}

\def\iff{\Leftrightarrow}

\def\<{\langle}

\def\>{\rangle}

\def\N{\mathbb{N}}

\def\Z{\mathbb{Z}}

\def\Q{\mathbb{Q}}

\def\F{\mathbb{F}}

\def\T{\mathfrak{T}}
\def\n{\mathfrak{n}}

\def\ns{\triangleleft}

\def\SL{\mathrm{SL}}

\def\isom{\cong}

\def\P{\mathbb{P}}

\def\Lamb{\Lambda'}
\def\flag{\mathcal{F}}

\def\L{\mathfrak{L}}

\newtheorem{theo}{Theorem}[section]

\newtheorem{lemma}[theo]{Lemma}

\newtheorem{cor}[theo]{Corollary}

\newtheorem{prop}[theo]{Proposition}

\theoremstyle{definition}

\newtheorem{Def}[theo]{Definition}

\newtheorem{ex}[theo]{Example}

\begin{document}

\title{The normal zeta function of the free class two nilpotent group
on four generators}

\author{Pirita Maria Paajanen
\\ Mathematical Institute\\24-29 St Giles'\\Oxford OX1 3LB\\UK\\
paajanen@maths.ox.ac.uk}

\maketitle

\begin{abstract}We calculate explicitly the normal zeta function
of the free group of class two on four generators, denoted by
$F_{2,4}$. This has Hirsch length ten.
\end{abstract}

\section{Introduction}

A zeta function of a group is a tool used in infinite group theory to understand subgroup
growth, the study of the properties of the sequence $(a_n(G))_{n\in\N}$, where $a_n(G)$
encodes the number of subgroups of index $n$ in a finitely generated infinite group $G$. By
studying analytic properties of the zeta function
\begin{equation*}\zeta^\leq_G(s)=\sum_{n=1}^\infty a_n^\leq(G)n^{-s},\end{equation*} where
$$a_n^\leq(G)=|\{H\leq G: |G:H|=n \}|,$$  we hope to get algebraic information about the
group.

We can define a normal zeta function $$\zeta^\ns_G(s)=\sum_{n=1}^\infty a_n^\ns(G)n^{-s}$$
 by the
additional condition that we count only normal
subgroups of $G$, so $$a_n^\ns (G)=|\{H\ns G: |G:H|=n \}|.$$

In the example in this paper we calculate a zeta function of a
finitely generated torsion-free nilpotent group --- we call these
the $\T$-groups --- and the nilpotency gives us the extra feature
that for $*\in \{\ns,\leq\}$ the zeta function decomposes as an
\emph{Euler product of local factors} (also known as the local
zeta functions) over primes $p$
$$\zeta_G^*(s)=\prod_p\zeta_{G,p}^*(s)$$ where
$$\zeta_{G,p}^*(s)=\sum_{n=0}^{\infty} a_{p^n}^*(G) p^{-ns}.$$

One of the main results in \cite{GSS} was to prove that these
local zeta functions are rational functions in $p^{-s}$. Moreover,
in \cite{GSS} it is proved that the zeta functions of the free
nilpotent groups of class two are uniform; this means that for
almost all, that is, for all but finitely many, primes $p$ the
local factor is the same for each prime. It should be noted that
this is not always the case, there are examples of zeta functions
whose local factors don't fall into even finitely many classes for
primes $p$. The elliptic curve example by du Sautoy discussed in
\cite{elliptic} and \cite{crelle} explains why the the local
factors of a zeta function depend on the points on an elliptic
curve modulo $p$.

More generally, in \cite{Annals} du Sautoy and Grunewald link the theory of zeta
functions of groups to the classical problem of counting points on varieties.

\begin{theo}[du Sautoy, Grunewald] Let $G$ be a finitely generated nilpotent infinite
group. Then there exist an algebraic variety $Y$ defined over
$\Q$, consisting of finitely many irreducible components $V_i$, $i
\in T$, defined over $\Q$. All these irreducible components are
smooth and intersect normally. There also exist rational functions
$P_I(x,y)\in\Q(x,y)$ for each $I\subseteq T$ with the property
that for almost all primes $p$ \begin{equation}
\zeta_{G,p}^*(s)=\sum_{I\subseteq T} d_{I,p}P_I(p,p^{-s}),
\end{equation} where
\begin{equation}d_{I,p}=|\{a\in \overline{Y}(\F_p):a\in\overline{V_i} \iff i\in I \}|
\end{equation} and $\overline{Y}$ means the reduction$\mod p$ of the variety $Y$.
 \end{theo}

However this theorem still has mainly theoretical meaning. One of the big questions
recently has been to describe which varieties appear in the decomposition. In \cite{elliptic}
and \cite{crelle} we saw that a zeta function of a certain group depends on the number of
$\F_p$-points on an elliptic curve. In this paper we can add Fano varieties of linear spaces to the list of
varieties we need to consider in this context. But even so, we need to be very careful, as
we shall see in this example. Not all points on the Fano variety of planes appear in the
zeta function. Only one of the two rulings of planes gets picked up by the zeta function.

In some sense this work can also be considered as a generalisation of Voll's work in
\cite{FE} and \cite{Amer}. In the paper \cite{FE} he proves that a local normal zeta
function satisfies a functional equation of the form $$\zeta_{G,p}^\ns(s)|_{p\mapsto
p^{-1}} =(-1)^{l} p^{m-ns} \cdot \zeta_{G,p}^\ns(s),$$ for some $l,m,n\in \N$, if the
Pfaffian hypersurface associated to the group is smooth, absolutely irreducible
 and does not contain lines.
The example that we shall present in this paper however involves a group whose associated Pfaffian is smooth, contains lines and
planes, and still satisfies a functional equation of the same form.

Apart from the pure interest in such a standard example as a nilpotent free group,
this calculation has application to the problem of counting $p$-groups. As du Sautoy
explains in \cite{IHES}, counting normal subgroups in the free class $c$ nilpotent
groups on $d$ generators is intimately related to counting $p$-groups of class $c$ on
$d$ generators --- the latter being identifiable by quotients of the former.

In addition to giving an explicit calculation of the normal zeta
function of $F_{2,4}$, we shall do it in a systematic and theoretical manner
which hopefully leads to a better theoretical understanding about the
shape of the function and its poles. In a future paper
\cite{smooth} I shall consider the more general case, where the
Pfaffian is assumed to be smooth.

This result was already announced in \cite{FE}.

It also should be noted that no computer calculations are used in this example. If one
tried to calculate this example using $p$-adic integrals, the cone conditions would
contain 53 polynomial equations and involve 55 variables. Thanks are due to Luke
Woodward who generated these conditions for me.

\subsection{Results}

Let $F_{2,d}$ be the free class two, $d$ generator nilpotent
group. $F_{2,d}= F_d/\gamma_3(F_d)$, where $F_d$ is the free group
on $d$ generators and $\gamma_3(F_d)$ is the third term in the
lower central series of $F_d$. $F_{2,d}$ has a presentation
$$F_{2,d}=\<x_1,\dots,x_d,y_1,\dots,y_{\frac{d(d-1)}{2}}:[x_1,x_2]=y_1,[x_1,x_3]=y_2,\dots,
[x_{d-1},x_d]=y_{\frac{d(d-1)}{2}}\>,$$ with the 
convention that generators not appearing
explicitly in the relations are assumed to commute.
Note that for this group $[F_{2,d},F_{2,d}]=Z(F_{2,d})$.

To start with, let us recall from \cite{GSS} what the normal zeta
functions of smaller examples of the class two nilpotent free
groups look like. The two generator one is isomorphic to the
Heisenberg group and has local normal zeta function
$$\zeta_{F_{2,2},p}^\ns(s)=\frac{1}{(1-p^{-s})(1-p^{1-s})(1-p^{2-3s})}.$$
The global zeta function $\zeta_{F_{2,2}}^\ns(s)=\prod_p \zeta_{F_{2,2},p}^\ns(s)$
has abscissa of convergence 2. The normal zeta
function of $F_{2,3}$ is
$$\zeta_{F_{2,3},p}^\ns(s)=\frac{1+p^{3-3s}+p^{4-3s}+p^{6-5s}+p^{7-5s}+p^{10-8s}}{(1-p^{-s})(1-p^{1-s})(1-p^{2-s})(1-p^{5-3s})(1-p^{8-5s})(1-p^{9-6s})}.$$
The global zeta function has abscissa of convergence 3.

In order to write the normal zeta function of $F_{2,4}$ in a
reasonable form ---  an explicit expression  would take two and
half pages of polynomial formula --- we need to introduce some
notation and terminology. We shall get back to these definitions
later and terminology will be made more precise.

\begin{Def}
The \emph{Igusa factor} with $n$ variables $\mathbf{U}=(U_1,\dots,U_n)$ is defined to be
$$I_n(\mathbf{U})=\sum_{I\subseteq \{1,\dots,n\}} b_I(p^{-1})\prod_{i\in
I}\frac{U_i}{1-U_i},$$ where
$b_I(p)$ is the number of $\F_p$-points on a flag variety of type $I$, this is a
polynomial in $p$, which can be written in terms of  $p$-binomial coefficients and is
independent of the prime $p$. \end{Def}

The precise definition of a flag variety of type $I$ will be given
in \ref{flag}.
\begin{ex}
In the above, we can write
$$\zeta_{F_{2,3},p}^\ns(s)=\zeta_{\Z^3,p}(s)\zeta_p(6s-9)
\cdot I_2(U_1,U_2),$$ where $U_1=p^{8-5s}$ and  $U_2=p^{5-3s},$ and
$\zeta_p(s)=\frac{1}{1-p^{-s}}$ denotes a local factor of the Riemann
zeta function.
\end{ex}

Let us denote by $d=h(F_{2,4}^{ab})=4$ the torsion-free rank of
the abelianisation and by $d'=h(Z(F_{2,4}))=6$ the torsion-free
rank of the centre, which is also an abelian group. In the general case of
free class two nilpotent groups $d$ will be number of generators
and $d'=\frac{d(d-1)}{2}$.

As will be explained later, to each class two nilpotent group we
can naturally associate a Pfaffian hypersurface in $\P^{d'-1}$.
Let $d_i$ be the dimension and $c_i$ the codimension of the Fano
variety of the $(i-1)$-dimensional linear subspaces on the
Pfaffian. Then $n_i=c_i+d_i$ is the dimension  of the variety of
$(i-1)$-dimensional linear subspaces in $\P^{d'-1}$.

\begin{Def}
We call a rational function of the form
\begin{equation*} E_i(X_i,Y_i) =
\frac{p^{-d_i}Y_i - p^{-n_i}X_i}{(1-X_i)(1-Y_i)}
\end{equation*} the $i^{th}$ \emph{exceptional factor} with variables $X_i$ and $Y_i$.
\end{Def}

With the above definitions we can now write down an expression for
the zeta function of $F_{2,4}$. It does show some structure, which
would be completely missing if the zeta function were to be
written explicitly.

\begin{theo}\label{main theorem}
The local normal zeta function of $F_{2,4}$ is for (almost) all
primes $p$
\begin{equation*}
\zeta_{F_{2,4},p}^\ns(s)=
\zeta_{\Z^4,p}(s)\zeta_p(10s-24)\cdot(W_0(p,p^{-s}) + \n_1(p)
W_1(p,
p^{-s})+\n_2(p)W_2(p,p^{-s})+\n_3(p)W_3(p,p^{-s}))\end{equation*}
where
\begin{align*}\n_1(p)&= (p^2 + 1)(p^2 + p + 1),\\ \n_{2}(p) & = (p + 1)(p^2 + 1)(p^2
+ p + 1),\\ \n_{3}(p)&= (p^2 + 1)(p + 1),
\end{align*}
and
$$W_i(X,Y)=I_{5-i}(X_{5},\dots,X_{i+1})E_i(X_i,Y_i)I_{i-1}(Y_{i-1},\dots,Y_1)$$
for $i=0,1,2,3$, with the convention
that $E_0=I_{-1}=I_0=1.$ %The parameter $t=2$ if $i=1$ and $t=1$ otherwise.
The numerical data are
\begin{align*}X_i&=p^{i(10-i)-(4+i)s}\\
Y_1&=p^{8-3s}\\ Y_2&=p^{13-5s}   \\ Y_3&=p^{15-6s}.
\end{align*}
\end{theo}

One can calculate the closed expression of this zeta function
where everything is put over a common denominator. However, it is
not done here, since the numerator is a polynomial with 490 terms
and the whole zeta function would take more than two pages to
write out. We can read off the degree of polynomial subgroup
growth from the expression above.

\begin{cor}
The abscissa of convergence of the global zeta function is $4$.
\end{cor}

This is still the torsion-free rank of the abelianisation.
However, for $F_{2,5}$ the abscissa is bigger than 5, which can be
deduced for instance from Theorem 1.3 in \cite{Abscissa}.

\begin{cor}
The function satisfies a functional equation of the form
$$\zeta_{F_{2,4},p}^\ns(s)|_{p\mapsto p^{-1}} =
p^{45-14s} \cdot \zeta_{F_{2,4},p}^\ns(s).$$
\end{cor}

This is of the form conjectured to hold for all class two normal
zeta functions, see section \ref{fe}.
\\
\\
\emph{Acknowledgements}. I would like to thank
Christopher Voll for numerous helpful conversations and
Marcus du Sautoy for his encouragement as a supervisor. 
My DPhil, part of which this paper forms, has been financially
supported by Academy of Finland, Finnish Academy
of Science and Letters, Osk. Huttunen Foundation and Oxford
University Graduate Scholarship Scheme.
I am also grateful for the friendliness of Helsinki University of
Technology during the writing up of this paper.

\section{Theory behind the calculation}

The class two nilpotent groups split naturally into two parts, the
derived group and the abelianisation. With free class two groups
we even have $G'=[G,G]=Z(G)$. By the Mal'cev correspondence there is a
Lie algebra $\L=\L(G)$ over $\Z$ constructed as the image under
$\log$ of $G$.

For a Lie algebra $\L$, and $*\in\{\leq,\ns\}$, we
define the zeta function $$\zeta_{\L}^*(s)=\sum_{n=1}^\infty
a_n^*(\L)n^{-s}$$ where $a_n^*(\L)$ is the number of subalgebras or
ideals of $\L$ of index $n$. Then we have the following theorem

\begin{theo}[\cite{GSS}]
Let $G$ be a finitely generated nilpotent group.
For $*\in \{\leq,\ns\}$ and almost all primes $p$
$$\zeta_{G,p}^*(s)=\zeta_{\L(G),p}^*(s).$$
\end{theo}

This allows us to concentrate on counting in a ring setting, which
is usually easier than the corresponding group setting.

According to \cite{GSS} (Lemma 6.1) for class two normal zeta
functions it is enough to run over the lattices of the centre of
the group (or Lie algebra) in order to calculate the zeta
function. In particular, if for each lattice $\Lambda'\subseteq
\L'$ we define $X(\Lambda')/\Lambda':=Z(\L/\Lambda'),$ then
\begin{align*}
\zeta_{\L}^\ns(s)=&\zeta_{\Z^d}(s)\sum_{\Lambda'\subseteq \L'}
|\L':\Lambda'|^{d-s}|\L:X(\Lambda')|^{-s}\end{align*}

We can also use the local version of this lemma by considering
$\L_p=\L\otimes_\Z \Z_p$. An ideal of a finite index in $\L_p$
corresponds to an ideal of a $p$-power index in $\L$. We shall
call a lattice $\Lamb\subseteq \Z_p^{d'}$ \emph{maximal} if
$p^{-1}\Lamb\not\subseteq \Z_p^{d'}$. It is enough to consider
only maximal lattices of $p$-power index, since if $\Lamb$ is not
maximal, then in general, $\Lamb= p^{r_{d'}}\Lamb_{max}$, where
$\Lamb_{max}$ is maximal in its class. Now
$$|\L'_p:\Lambda'|=p^{d'r_{d'}} |\L'_p:\Lambda'_{max}|$$ and
$$|\L_p:X(\Lambda')|=p^{dr_{d'}}|\L_p:X(\Lambda'_{max})|.$$

Then we have
\begin{equation*}\begin{split}
\zeta_{\L,p}^\ns(s)=
\zeta_{\L_p}^\ns(s)=&\zeta_{\Z_p^d}(s)\sum_{\Lambda'\subseteq
\L_p'} |\L_p':\Lambda'|^{d-s}|\L_p:X(\Lambda')|^{-s}\\ =&
\zeta_{\Z_p^d}(s)\zeta_p((d+d')s-dd')A(p,p^{-s}),\end{split}
\end{equation*} where 
\begin{equation}\label{generating functions}
A(p,p^{-s})=\sum_{\underset{\Lamb
\text{maximal}}{ \Lambda'\subseteq
\L_p'}}|\L_p':\Lambda'|^{d-s}|\L_p:X(\Lambda')|^{-s}.\end{equation}

Recall from \cite{Amer} the definition of the weight functions
\begin{align*}w(\Lambda')&:=log_p(|\L'_p:\Lambda'|)\\
w'(\Lambda')&:=w(\Lamb)+\log_p(|\L_p:X(\Lamb)|).\end{align*}

Put $T:=p^{-s}$. Then we can write our generating function as
\begin{equation}\label{generating function}A(p,T)=\sum_{\Lambda'\subseteq \L'_p, maximal}
p^{dw(\Lambda')}
T^{w'(\Lambda')}.\end{equation}

Let us consider $\L_p$ and a lattice $\Lamb\in \L_p'\isom
\Z_p^{d'}.$ In order to do the counting effectively, we enumerate
the lattices in the centre using the elementary divisor type. The
lattices maximal in their class have elementary divisor types
$(1,p^{r_{d'-1}},p^{r_{d'-2}+r_{d'-1}},\dots,p^{r_1+r_2+\dots+r_{d'-1}})$,
where $r_i\geq 0$. Let us denote by $\{i_1,\dots,i_n\}_<$ the
ordered set where $i_1<i_2<\dots i_l$. Then the type of lattice is
given by $\nu(\Lambda')=(I,r_I)$ where $$I=\{i_1,\dots
i_l\}_{<}\subseteq\{1,\dots,d'-1\}, r_I=(r_{i_1},\dots,r_{i_l})$$
and $I=\{i_j\in\{1,\dots,d'-1\}:r_{i_j}>0\}$.

This way we have the same enumeration for the lattices, and for
the weight functions depending on the lattices' types. From these
definitions and the above notation for the lattices we have
immediately $$w(\Lamb)=\sum_{i\in I} ir_i.$$

\begin{Def}\label{flag}
A \emph{flag} of type $I$ in $\P^{d'-1}(\F_p)$,
$I\subseteq\{1,\dots,d'-1\}$ is a sequence $(V_{i_j})_{i_j\in I}$
of incident vector spaces
$$\P^{d'-1}(\F_p)>V_{i_1}>\dots>V_{i_l}>\{0\}$$ with
$\mathrm{codim}_{\F_p}(V_{i_j})=i_j-i_{j-1}$ in the vector space
$V_{i_{j-1}}$ . The flags of type $I$ form a projective variety
$\flag_I$, whose number of $\F_p$-points is given by
$b_I(p)\in\Z[p]$ a polynomial whose leading term equals
$p^{\mathrm{dim} \flag_I}$ and can be easily expressed as a
product of $p-$binomial coefficients. Using duality properties of
$p$-binomial coefficients we observe that
$$b_I(p^{-1})=p^{-\dim\flag_I} b_I(p).$$\end{Def}

$\SL_{d'}(\Z_p)$ acts transitively on maximal lattices of an
elementary divisor type
$$(1,p^{r_{d'-1}},p^{r_{d'-2}+r_{d'-1}},\dots,p^{r_1+r_2+\dots+r_{d'-1}})$$
relative to the standard lattice.

Let $G_\nu$ denote the stabiliser of the diagonal matrix
$diag(1,p^{r_{d'-1}},p^{r_{d'-2}+r_{d'-1}},\dots,p^{r_1+r_2+\dots+r_{d'-1}})$
in $\SL_{d'}(\Z_p)$.  Then by the orbit-stabiliser theorem we have
$$\{\text{maximal lattices of type}\
\nu\}\overset{1-1}{\leftrightarrow} \SL_{d'}(\Z_p)/G_\nu.$$

Thus we identify a maximal lattice $\Lamb$ with the pair
$(\nu,\alpha G_\nu)$ consisting of its partition of elementary
divisors and the respective coset of the stabiliser $G_\nu.$

Let $\alpha\in\SL_{d'}(\Z_p)$ and $I$ be a type of a lattice
$\Lambda'$. Let $\overline{\alpha}$ denote the reduction $\mod p$
and define vector spaces for $i_j\in I$, $I=\{i_1,\dots,i_l\}_<,$
$$V_{i_j}=\<\overline{\alpha^{i_1}},\dots,\overline{\alpha^{i_j}}  \>
<\P^{d'-1}(\F_p)$$ We observe that $\mathrm{codim}_{\F_p}(V_{i_j})=i_j-i_{j-1}.$

Following Voll we call the flag $(V_i)_{i\in I}$ of type $I$ the
flag associated with $\Lambda'$ if $\nu(\Lambda')=I$ and
$\Lambda'$ corresponds to the coset $\alpha G_\nu$ under the 1-1
correspondence.

\begin{Def}\label{lift}
Given a fixed flag variety $\mathcal{F}_I$ of type $I=\{i_1,\dots,i_l\}_<$  on $\P^{d'-1}$, then we
call a maximal lattice $\Lambda'$ of type $J=\{j_1,\dots,j_k\}_<$ where $k\geq l$  a \emph{lift}
of
$\mathcal{F}_I$
if its associated flag $(V_{j})_{j\in J}$
contains $\mathcal{F}_I$, $I\subseteq J$ and $i_m=j_m$ for $m=1,\dots,l.$\end{Def}

The reason we define these is that we can describe $|\L_p:X(\Lamb)|$ by set of
polynomial equations.

\begin{theo}[Voll \cite{FE}]
Let $\Lambda'$ correspond to the coset $\alpha G_\nu$ under the correspondence defined above,
where $\alpha\in \SL_6(\Z_p)$ with column vectors $\alpha^j$, $j=1,\dots,6$.
Then $|\L_p:X(\Lambda')|$
equals the index of the kernel of the following system of linear congruences in
$\L_p/\L'_p$:
\begin{equation}\label{congruences}
\forall i\in\{1,\dots,6\}\ \ \bar{\mathbf{g}}M(\alpha^j)\equiv 0 \mod p^{r_i+\dots+r_6}
\end{equation} where $\bar{\mathbf{g}}=(\bar{g_1},\dots,\bar{g_4})\in \L_p/\L'_p
\isom \Z_p^4.$
\end{theo}

From this it is clear that we need to separate cases when
$M(\alpha^j)$ is singular and when it is not. The space of
singular matrices can be described using the Pfaffian
hypersurface, which is defined by the vanishing of the determinant
of the matrix of relations $\det(M(\mathbf{y}))=0$.

In addition to the above two weight functions we need to define
the multiplicity function $\mu$ which measures how many lattices
there are of each given type.

\begin{Def} Let $\Lamb$ be a lattice of type $I$.
The multiplicity of $\Lambda'$, denoted by $\mu(\Lambda')$ is the
number of lattices of type $I$ divided by the number of points on
the flag variety of type $I$.\end{Def}

We can calculate the multiplicity $\mu(\Lamb)$ by defining the
function mu (as a `measure´ or a multiplicity) to measure the size
of the set of $x\in p\Z_p/(p^{a})$ of a fixed $p$-adic valuation
as follows:

\begin{Def}
Let $a,b$ be positive integers. We define a function $\mu$ of $a$
and $b$ as \begin{equation} \mu(a,b)=|\left\{x\in
p\Z_p/(p^{a}):v_p(x)=b\right\}|=\begin{cases} 1 & \text{if $a=b$}\\
p^{a-b}(1-p^{-1}) & \text{if $a>b$}\\ 0
&\text{otherwise.}\end{cases}
\end{equation}
This definition extends naturally also to a vector $\mathbf{b}=(b_1,\dots,b_n)$.
We denote $\mu(a,{\bf b})=|\left\{\mathbf{x}\in
p\Z_p/(p^{a}):v_p(x_i)=b_i\right\}|=
\mu(a;b_1,b_2,b_3,\dots,b_n)=\mu(a,b_1)\mu(a,b_2)\dots\mu(a,b_n).$
\end{Def}

Now we can calculate the multiplicity of lattices of a given type $I$.

\begin{lemma}\label{no of lattices}
Let $\Lamb$ be a lattice of type $I$  corresponding to the coset
$\alpha G_\nu$ under the 1-1 correspondence. Let $a_{ij}$ be the
$ij$-entry of $\alpha$ (in the ith column from the right and
(j-1)th row from the bottom). Then the multiplicity of $\Lamb$ is
$$\mu(\Lambda') =\prod_{\underset{j\leq i}{\underset{j\in I}{i\in
I}}} \sum_{a_{ij}=1}^{\sum_{k=i}^j r_k }\mu (\sum_{k=i}^j r_k,
a_{ij})= p^{-\dim\flag_I}p^{\sum_{i\in I} (d'-i)ir_i}.$$
\end{lemma}

\begin{proof}
A simple computation using the orbit-stabiliser theorem.
\end{proof}

\section{The group and its Pfaffian hypersurface}

The free class two nilpotent group on 4 generators has the
following presentation \begin{equation*}\begin{split}
F_{2,4}=&\<x_1,x_2,x_3,x_4,y_1,y_2,y_3,y_4,y_5,y_6:
[x_1,x_2]=y_1,[x_1,x_3]=y_2,\\&[x_1,x_4]=y_3,
[x_2,x_3]=y_4,[x_2,x_4]=y_5,[x_3,x_4]=y_6\>.\end{split}\end{equation*}
We can also write the commutator structure in the form of the
matrix of relations as

$$M(\bf{y})=\begin{pmatrix}0 & y_1 &y_2&y_3\\ -y_1&0 &y_4&y_5\\
-y_2&-y_4&0&y_6\\ -y_3&-y_5&-y_6&0
\end{pmatrix}.$$

Then the presentation is given by
 \begin{equation*}\begin{split}
F_{2,4}=&\<x_1,x_2,x_3,x_4,y_1,y_2,y_3,y_4,y_5,y_6:
[x_i,x_j]=M(\mathbf{y})_{ij}\>.\end{split}\end{equation*}

The Pfaffian hypersurface is defined by the square-root of the  determinant of the matrix of relations
$M(\mathbf{y})$; in this case it is defined by the equation
\begin{equation}y_1y_6-y_2y_5+y_3y_4=0.\end{equation}
This is a smooth quadric four-fold in $\P^5$. It is a fact (see
e.g. Theorem 22.13 in \cite{Harris}) that a smooth quadric of
dimension $m$ contains linear subspaces only up to dimension less
than or equal to $m/2$, thus the linear subspaces lying on the
Pfaffian hypersurface $\mathfrak{P}$ associated with $F_{2,4}$ are
points, lines and planes, these are also smooth by the above
mentioned Theorem 22.13 . As we can find, for example from  from
\cite{Hirschfeld} p. 5, there are $\mathfrak{P}_1(p)= (p^2 +
1)(p^2 + p + 1)$ points, $ \mathfrak{P}_{2}(p)  = (p + 1)(p^2 +
1)(p^2 + p + 1)$ lines and  $\mathfrak{P}_{3}(p)= 2(p^2 + 1)(p +
1)$ planes over $\F_p$. There are two rulings of planes, which
again is from the general theory of quadrics, see Theorem 22.14 in
\cite{Harris}.

Thus the different lifts of fixed flags that lie on the Pfaffian
hypersurface (Recall Definition \ref{lift}) and can cause
different congruence conditions and weight functions are the
following:

\begin{enumerate}
\item The lattice $\Lamb$ lifts the flag
$\<\alpha^1\>$, which consists of a point on the Pfaffian. Note
that the line $\<\alpha^1,\alpha^2\>$ does not lie on the quadric
four-fold in this case.

\item \begin{itemize} \item[(a)]  The lattice $\Lamb$ lifts the flag
$\<\alpha^1,\alpha^2\>$, which
consists of a line  on the Pfaffian.
   \item[(b)] $\Lamb$ lifts the flag  $\<\alpha^1,\alpha^2\>\geq \<\alpha^1\> $, which is a
line-point pair on the Pfaffian. \end{itemize}

\item \begin{itemize} \item[(a)] $\Lamb$ lifts the flag $\<\alpha^1,\alpha^2,\alpha^3\>$, which
consists of a plane
on the Pfaffian.
\item[(b)] $\Lamb$ lifts the flag  $\<\alpha^1,\alpha^2,\alpha^3\>\geq \<\alpha^1,\alpha^2\>$,
 which
consists of a plane-line pair on the Pfaffian.
\item[(c)] $\Lamb$ lifts the flag $\<\alpha^1,\alpha^2,\alpha^3\>\geq \<\alpha^1\>$, consisting of a
plane-point pair on the Pfaffian.
\item[(d)] $\Lamb$ lifts the flag $\<\alpha^1,\alpha^2,\alpha^3\>\geq \<\alpha^1,\alpha^2\>\geq
\<\alpha^1\>$, which is a fixed plane-line-point triplet on the
Pfaffian.
\end{itemize}
\end{enumerate}

\subsection{Weight functions}

In this section we consider different types of lattices and what kind of
weight functions $w'(\Lamb)$ we
obtain following the above geometric list.

The easiest case is when a lattice lifts a flag no part of which
lies on the Pfaffian hypersurface, so
all the matrices $M(\alpha^i)$ are non-singular and the congruence conditions reduce to:

$$\forall i\in\{1,\dots,6\}\ \ \bar{\mathbf{g}}\equiv 0 \mod p^{r_i+\dots+r_5}.$$

The weight function $w'(\Lamb)=\sum_{i\in I} (4+i)r_i$ depends
only on the lattice's type.

Now we need to consider the conditions on the Pfaffian as listed
above. The easiest of them is that a lattice $\Lamb$ lifts a fixed
point $x\in \P^5(\F_p)$ on the Pfaffian hypersurface.

\begin{lemma}[Voll, Proposition 4 in \cite{FE}]
Let $x\in \mathfrak{P}$, and let $\Lamb$ be a lattice lifting this fixed point. Then
the weight function
is $w'(\Lamb)=\sum_{i\in I} (4+i)r_i-2\min\{r_1,v_p(a_{11})\}$.
This depends on more than the lattice's
type, but it is independent of the point $x$ chosen.
\end{lemma}

\begin{proof}(Proposition 4 in \cite{FE})

Let $\alpha^1=(a_{11},a_{12},\dots,a_{16})\in \P^{5}(\Z_p/p^{r_1})$ and
$\alpha^1=(0,1,\dots,1) \mod p.$
Then we can choose
local coordinates such that around any of the $\n_1(p)$ points of the Pfaffian hypersurface
$\mod p$ the congruence conditions look like
\begin{align*}\overline{\mathbf{g}}\cdot \begin{array}({c c c c })0 &a_{11} &0&0\\
-a_{11}&0&0&0\\ 0&0&0&1\\0&0&-1&0 \end{array}
&\equiv 0 {\mod p^{r_1+\dots+r_5}}\\
\overline{\mathbf{g}}&\equiv 0 {\mod p^{r_2+\dots+r_5}}.
\end{align*}

It can be read off that the weight function is
$w'(\Lamb)=\sum_{i\in I} (4+i)r_i-2\min\{r_1,v_p(a_{11})\}$.

\end{proof}

With lines we have the following situation: The variety of lines
is smooth and irreducible, and all lines belong to the same
family, so it is enough to consider one given line only.

\begin{lemma} The weight function for lattices lifting a given line on the Pfaffian is
$w'(\Lamb)=6r_2-\min\{r_2,v_p(a_{15}),v_p(a_{25}),v_p(a_{24})\}$. This is independent of
the line chosen. \end{lemma}

\begin{proof}
The lattice of type $(1,1,1,1,p^{r_2},p^{r_2})$ is in one-to-one
correspondence with the pair of vectors
$\alpha^1=(a_{15},a_{14},a_{13},a_{12},0,1)^t$
and $\alpha^2=(a_{25}, a_{24},a_{23},a_{22},1,0)^t$,
where $a_{ij}\in p\Z_p/(p^{r_2})$, so that $\mod p $ the span
$\<\alpha^1,\alpha^2\>$ defines a line on the Pfaffian.
With suitable local coordinate changes the congruences reduce to

\begin{equation*}
\overline{\mathbf{g}}
\begin{array}({c c c c}) 0 & a_{15} &0  & 0 \\
0 & 0 & 0 & 0
\\0  &0&0&1\\
0 &0&-1&0\\
\end{array}\equiv 0 \pmod{p^{r_2}}
\end{equation*}

\begin{equation*}
\overline{\mathbf{g}}
\begin{array}({c c c c}) 0 & a_{25} & a_{24}  & 0 \\
0 & 0 & 0 & 1
\\0  &0&0&0\\
0 &0&0&0\\
\end{array}\equiv 0 \pmod{p^{r_2}}.
\end{equation*}

\end{proof}

\begin{cor}
For a mixed lattice of type $(1,\dots,1,p^{r_2},p^{r_1+r_2})$ any
flag $\<\alpha^1\>\leq\<\alpha^1,\alpha^2\>$ gives the same weight
function.
\end{cor}

\begin{proof}
As each of the lines and each of the points give the same weight
function, so does any flag of any point-line combination. We can
thus take the following vectors
$\alpha^1=(a_{15},a_{14},a_{13},a_{12},a_{11},1)^t$ and
$\alpha^2=(a_{25}, a_{24},a_{23},a_{22},1,0)^t$, where $a_{2j}\in
p\Z_p/(p^{r_2})$, $a_{1j}\in p\Z_p/(p^{r_1+r_2})$ for $2\leq j\leq
5$, $a_{11}\in p\Z_p/(p^{r_1})$. With suitable chance of local
coordinates the congruence conditions look like

\begin{equation*}
\overline{\mathbf{g}}
\begin{array}({c c c c}) 0 & a_{15} &0  & 0 \\
-a_{15} & 0 & 0 & 0
\\0  &0&0&1\\
0 &0&-1&0\\
\end{array}\equiv 0 \pmod{p^{r_1+r_2}}
\end{equation*}

\begin{equation*}
\overline{\mathbf{g}}
\begin{array}({c c c c}) 0 & a_{25} & a_{24}  & 0 \\
0 & 0 & 0 & 1
\\0  &0&0&0\\
0 &0&0&0\\
\end{array}\equiv 0 \pmod{p^{r_2}},
\end{equation*}
from which we can read off the weight function to be
$w'(\Lamb)=6r_2+5r_1-\min\{r_1,v_p(a_{15})\}-\min\{r_1+r_2,v_p(a_{15}),v_p(a_{25})+r_1,v_p(a_{24})+r_1\}.$
\end{proof}

However, with planes we need to be more careful. As planes are the
highest dimensional linear subspace on the quadric, there are two
families of planes, and as we shall see, these rulings do not
behave equally; only one of them gives a different weight
function. The variety of the planes is smooth, thus it is enough
to consider a representative in each one of the rulings.

%of both rulings, other planes in the same ruling will behave uniformly. %\footnote{Terminology
%and claims here can be complete %rubbish. I need to check these things properly.}

\begin{lemma}
The plane generated by the vectors
$\alpha^1=(1,0,0,a_{33},a_{34},a_{35})^t$,
$\alpha^2=(0,1,0,a_{23},a_{24},a_{25})^t$ and
$\alpha^3=(0,0,1,a_{13},a_{14},a_{15})^t$, $a_{ij}\in
p\Z_p/(p^{r_3})$  lies on the Pfaffian $\mod p$, but doesn't give
a different weight function. For this half of the planes the
weight function is $w'(\Lamb)=7r_3$ which is also the weight
function for planes outside the Pfaffian. There are $(p^2+1)(p+1)$
such planes.
\end{lemma}

\begin{proof}
It is easy to see that the plane spanned by $\<\alpha^1,\alpha^2,\alpha^3\>$ is on the Pfaffian $\mod p$. The congruence conditions $\mod p$ are $\overline{\mathbf{g}}\cdot M(\alpha^i)\equiv 0\mod p$ for $i=1,2,3,$ and

\begin{equation*}
\overline{\mathbf{g}}
\begin{array}({c c c c}) 0 & 1 &0  & 0 \\
-1 & 0 & a_{33} &a_{34}
\\0  &-a_{33}&0&a_{35}\\
0 &-a_{34}&-a_{35}&0\\
\end{array}\equiv 0 \pmod{p^{r_3}}
\end{equation*}

\begin{equation*}
\overline{\mathbf{g}}
\begin{array}({c c c c}) 0 & 0 &1  & 0 \\
0 & 0 & a_{23} &a_{24}
\\-1  &-a_{23}&0&a_{25}\\
0 &-a_{24}&-a_{25}&0\\
\end{array}\equiv 0 \pmod{p^{r_3}}
\end{equation*}

\begin{equation*}
\overline{\mathbf{g}}
\begin{array}({c c c c}) 0 &0 &0 &1 \\
0 & 0 & a_{13} &a_{14}
\\0&-a_{13}  &0&a_{15}\\
-1 &-a_{14}&-a_{15}&0\\
\end{array}\equiv 0 \pmod{p^{r_3}}.
\end{equation*}

This set of equations has rank 4, and thus the weight function is $w'(\Lamb)=7r_3$ which is also the weight
function for planes outside the Pfaffian. \end{proof}

\begin{lemma} For the other $(p^2+1)(p+1)$ of planes the weight function is
$w'(\Lamb)=7r_3-\min\{r_3,v_p(a_{15}),v_p(a_{25}),v_p(a_{24}),v_p(a_{35}),v_p(a_{34}),v_p(a_{33})\}.$
\end{lemma}

\begin{proof}
The other family of planes can be represented by the span of the three vectors
$\alpha^1=(a_{15},a_{14},a_{13},0,0,1)^t$, $\alpha^2=(a_{25}, a_{24},a_{23},0,1,0)^t$,
$\alpha^3=(a_{35},a_{34},a_{33},1,0,0),$ where $a_{ij}\in p\Z_p/(p^{r_3})$. With a suitable change of local
coordinates the congruence conditions look like.

\begin{equation*}
\overline{\mathbf{g}}
\begin{array}({c c c c}) 0 & a_{15} &0  & 0 \\
0 & 0 & 0 & 0
\\0  &0&0&1\\
0 &0&-1&0\\
\end{array}\equiv 0 \pmod{p^{r_3}}
\end{equation*}

\begin{equation*}
\overline{\mathbf{g}}
\begin{array}({c c c c}) 0 & a_{25} & a_{24}  & 0 \\
0 & 0 & 0 & 1
\\0  &0&0&0\\
0 &0&0&0\\
\end{array}\equiv 0 \pmod{p^{r_3}}
\end{equation*}

\begin{equation*}
\overline{\mathbf{g}}
\begin{array}({c c c c}) 0 & a_{35} &a_{34}  &a_{33} \\
0 & 0 & 0 & 0
\\0  &0&0&0\\
0 &0&0&0\\
\end{array}\equiv 0 \pmod{p^{r_3}}
\end{equation*}
We can read off the weight function to be as claimed.
\end{proof}

\begin{cor}
For the mixed lattices containing the type $(1,\dots,p^{r_3},p^{r_2+r_3},p^{r_1+r_2+r_3})$
$r_3>0$ and $r_1,r_2\geq 0$ we need to consider the planes that give a different weight
function only. Other planes reduce this to the two-dimensional case. \end{cor}

\begin{proof}
For mixed lattices corresponding to the vectors
$\alpha^3=(1,a_{34},a_{33},a_{32},a_{31},a_{30})$, $\alpha^2=(0,1,a_{23},a_{22},a_{21},a_{20})$,
$\alpha^1=(0,0,1,a_{12},a_{11},a_{10})$
the rank of matrices on the level of planes is four.
\end{proof}

\begin{lemma}
The weight function for mixed lattices where the plane is not one of the rank four planes is
\begin{equation}\begin{split}
w'(\Lamb)=&7r_3+6r_2+5r_1-\min\{r_1,v_p(a_{15})\}\\ -\min
\{&r_1+r_2+r_3, v_p(a_{15}),
v_p(a_{25})+r_1,v_p(a_{24})+r_1,\\
&v_p(a_{35})+r_1+r_2,v_p(a_{34})+r_1+r_2,
v_p(a_{33})+r_1+r_2\}.\end{split}\end{equation}
\end{lemma}

\begin{proof}
We count over lattices of elementary divisor type
$(1,1,1,p^{r_3},p^{r_2+r_3},p^{r_1+r_2+r_3})$. The lattices
$\Lamb$ lifting this flag are in one-to-one correspondence with
the three of vectors encoded as columns of the matrix. Call the
vectors $\alpha^3$, $\alpha^2$ and $\alpha^1$, respectively,\\

$\begin{pmatrix}a_{35} & a_{25} &a_{15}\\ a_{34}&a_{24}&a_{14}\\ a_{33}
&a_{23}&a_{13}\\ 1 &
a_{22} & a_{12}\\ 0 &1&a_{11}\\ 0& 0&1
\end{pmatrix}$\\

where $a_{35},a_{34},a_{33}\in p\Z_p/(p^{r_3})$, $a_{25},a_{24},a_{23}\in
p\Z_p/(p^{r_2+r_3})$, $a_{22}\in p\Z_p/(p^{r_2})$, $a_{15},a_{14},a_{13}\in
p\Z_p/(p^{r_1+r_2+r_3})$, $a_{12}\in p\Z_p/(p^{r_1+r_2})$, $a_{11}\in
p\Z_p/(p^{r_1}).$

So what we are required to do is to solve the following
congruences:

\begin{equation}\begin{split}
\overline{\mathbf{g}} \cdot M(\alpha^3)&\equiv 0 \pmod{p^{r_3}}\\\overline{\mathbf{g}}\cdot M(\alpha^2)&\equiv 0
\pmod{p^{r_2+r_3}}\\ \overline{\mathbf{g}}\cdot M(\alpha^1)&\equiv 0
\pmod{p^{r_1+r_2+r_3}}\\
\end{split}\end{equation}
simultaneously.

With suitable local coordinate changes the conditions reduce to

\begin{align*}
\overline{\mathbf{g}}
\begin{array}({c c c c}) 0 & a_{15} &0  & 0 \\
-a_{15} & 0 & 0 & 0
\\0  &0&0&1\\
0 &0&-1&0\\
\end{array}&\equiv 0 \pmod{p^{r_1+r_2+r_3}}\\
\overline{\mathbf{g}}
\begin{array}({c c c c}) 0& a_{25} & a_{24} & 0
\\ 0 &0 & 0 & 1
\\0&0 &0 &0\\
0  & 0& 0& 0\\
\end{array}&\equiv 0 \pmod{p^{r_2+r_3}}\\
\overline{\mathbf{g}}
\begin{array}({ c c c c }) 0& a_{35} & a_{34} & a_{33}
\\0& 0& 0&0\\
0  & 0& 0& 0\\0&0&0&0
\end{array}&\equiv 0 \pmod{p^{r_3}}
\end{align*}

The explicit weight function for this is
\begin{equation}\begin{split}
w'(\Lamb)=&7r_3+6r_2+5r_1-\min\{r_1,v_p(a_{15})\}\\ -\min
\{&r_1+r_2+r_3, v_p(a_{15}),
v_p(a_{25})+r_1,v_p(a_{24})+r_1,\\
&v_p(a_{35})+r_1+r_2,v_p(a_{34})+r_1+r_2,
v_p(a_{33})+r_1+r_2\}.\end{split}\end{equation}
\end{proof}

Other types of mixed lattices containing $r_3$ are obtained by
putting either of $r_1$ or $r_2$ equal to zero in the formula
above.

\section{Generating functions}

Recall from (\ref{generating functions}) that the zeta function
takes the shape
\begin{equation}\zeta_{F_{2,4}}^\ns(s)=\zeta_{\Z^4,p}(s)\zeta_p(10s-24)\cdot A(p,p^{-s}).
\end{equation}

We can decompose the generating function $A(p,p^{-s})$ further to run over
lattices of fixed type, and write it as \begin{equation}\label{generating
}A(p,p^{-s})=\sum_{I\subseteq\{1,\dots,d'-1\}}A_I(p,p^{-s})\end{equation} where
\begin{equation}\label{actual
function}A_I(p,p^{-s})=\sum_{\nu(\Lamb)=I}p^{dw(\Lamb)-sw'(\Lamb)}\mu(\Lamb),\end{equation}
$\Lamb$ is a representative lattice of type $I$, and $\mu(\Lamb)$
is its multiplicity.

However, the above is still not the optimal way to decompose
things; we somehow need to involve the dependence on the geometric
pieces which give a different weight function. A more subtle
decomposition is needed.

Let us denote by $\{i^*\}$, $i=1,2,3$  those lattices of type
$\{i\}$  that lift a given $i-1$-dimensional linear subspace of
the projective space $\P^5$ that lies on the Pfaffian and
accordingly by $A_{i^*}$ the generating function counting along
these lattices. Similarly for mixed lattices. For instance, we
have the lattices $\{1^*,2^*\}$ which lift a flag of type
$\flag_{1,2}$ consisting of a line containing a point on the
Pfaffian, and $A_{1^*,2^*}$ for the corresponding generating
function. In short, $A_I$ denotes the generating function along
lattices of type $I$ and the extra stars on elements $I$ just tell
if this part of the lattice has lifted a flag variety that lies on
the Pfaffian hypersurface.

\subsection{Indexing}

The indexing in this case has to be done carefully in order to take all the possibilities
into account and we shall use the indexing set $I\subset\{1,2,3,4,5,1^*,2^*,3^*\}$, such
that if $i\in I$ then $i^*\not\in I$ and if $j^*\in I$ then $k\not\in I$ when $k\leq j.$

Now the generating function $A(p,p^{-s})$ can be split into parts
\begin{align}
A(p,p^{-s})&=\sum_{I\subseteq\{1,\dots,5\}}c_{I,p}A_I(p,p^{-s})+
\sum_{\underset{J_2\subseteq\{2,\dots,5\}}{I=1^*\cup
J_2}}c_{I,p}A_I(p,p^{-s})
+\sum_{\underset{J_2\subseteq\{3,4,5\}}{\underset{J_1\subseteq\{1^*\}}{I=J_1\cup
2^*\cup J_2}}}c_{I,p}A_I(p,p^{-s})\nonumber \\&
+\sum_{\underset{J_2\subseteq\{4,5\}}{\underset{J_1\subseteq\{1^*,2^*\}}
{I=J_1\cup 3^*\cup J_2}}}c_{I,p}A_I(p,p^{-s}),\label{decomposition}
\end{align} where $A_I(p,p^{-s})$ are as in (\ref{actual function}) and $c_{I,p}$ are
coefficients depending on the type of lattice and number of $\F_p$-points on certain
varieties, and can be explicitly written down as polynomials in $p$.

In order to write down the coefficients $c_{I,p}$ we need some definitions. Let $\n_i(p)$
be the number of $\F_p$-rational points of the Fano varieties of $(i-1)$-dimensional
subspaces on the Pfaffian hypersurface. In the case of planes, however, we need to take the
number of planes on the ruling that gave the different weight function. Explicitly
\begin{align*}
\n_1(p)&= (p^2 + 1)(p^2 + p + 1),\\ \n_{2}(p) & = (p + 1)(p^2 +
1)(p^2 + p + 1),\\ \n_{3}(p)&= (p^2 + 1)(p + 1).
\end{align*}
So in this example we have $\n_3(p)=(p^2+1)(p+1)$ but really there are
$\mathfrak{P}_3(p)=2(p^2+1)(p+1)$ $\F_p$-points on the Fano variety of planes on a quadric
four-fold.
%Apart from the planes, all other Fano varieties behave uniformly.

Now we can explicitly describe the coefficients $c_{I,p}$
appearing in the generating function. Let us denote by $b_I(p)$
the number of points on the flag variety defined by lattices of
type $I$. We also write $I- k$ to mean the type of lattice we get
if we subtract from each index  $i\in I$, the number $k$, so if
$I=\{4,5,6\}$ then $I- 3=\{1,2,3\}$. With this notation we have
\\

1) If $I\subseteq \{1,\dots,5\}$ and $I=\{i_1,\dots,i_n\}$, then
\begin{align*} c_{I,p}&=
\binom{6-i_1-\dots-i_{n-1}}{i_{n-1}-i_n}_p\dots
\binom{6-i_1}{i_1-i_2}_p\left(\binom{6}{i_1}_p-\n_{i_1}(p)\right)\\&
=b_{I}(p)-b_{I- i_1}(p)\n_{i_1}(p).
\end{align*}

2) If $I=\{i_1,\dots,i_n,k^*,j_1,\dots,j_r\}$, then
\begin{align*}c_{I,p}&=\binom{6-k-\dots-j_{r-1}}{j_r-j_{r-1}}_p\dots
\binom{6-k-j_1}{j_2-j_1}_p\cdot\\&\cdot
\left(\binom{6-k}{k-j_1}_p\n_k(p)\binom{k-i_1\dots -i_n
}{i_n-i_{n-1}}_p\dots\binom{k}{i_1}_p-\n_{j_1}(p)\binom{j_1-i_1-\dots
i_n }{k-i_n}_p \dots\binom{j_1}{i_1}_p \right)\\&= b_{J_2-
k}(p)\n_k(p)b_{J_1}(p)-b_{J_2-(k-j_1)}(p)\n_{k+1}(p)b_{J_1\cup
k}(p).
\end{align*}

To see exactly where in the zeta function we have the dependence
on the $\F_p$-points of a Fano variety, let us now rearrange
(\ref{decomposition}) into the form
\begin{equation}A(p,p^{-s})=W_0(p,p^{-s})+\n_1(p)W_1(p,p^{-s})+\n_2(p)W_2(p,p^{-s})
+\n_3(p)W_3(p,p^{-s}).\end{equation}

From the formulae for the coefficients $c_{I,p}$ we see that the
only $c_{I,p}$ that don't depend on any of the $\n_i(p)$ come from
the first summand, and so we get

\begin{equation}
W_0(p,p^{-s})=\sum_{I\subseteq\{1,\dots,5\}}
b_I(p)A_I(p,p^{-s}).\end{equation}

We also observe that $\n_1(p)$ appears only in the first two
summands in (\ref{decomposition}) and thus

\begin{align*}
W_1(p,p^{-s})&=\sum_{\underset{J_2\subseteq\{2,\dots,5\}}{I=1^*\cup
J_2}}b_{I- 1}(p)A_I(p,p^{-s}) -
\sum_{\underset{J_2\subseteq\{2,\dots,5\}}{I=1\cup
J_2}}b_{I-1}(p)A_I(p,p^{-s}) \\&
=\sum_{I\subseteq\{2,\dots,5\}}b_{I- 1}(p)(A_{1^*\cup
I}(p,p^{-s})-A_{1\cup I}(p,p^{-s})).
\end{align*}
Similarly we extract the dependence on $\n_2(p)$ and notice that
it appears only in the first two summands of (\ref{decomposition})

\begin{align*}
W_2(p,p^{-s})&=\sum_{\underset{J_2\subseteq\{3,4,5\}}{\underset{J_1\subseteq\{1^*\}}{I=J_1\cup
2^*\cup J_2}}}b_{J_2-2}(p) b_{J_1}(p) A_I(p,p^{-s})-
\sum_{\underset{J_2\subseteq\{3,4,5\}}{\underset{J_1\subseteq\{1^*\}}{I=J_1\cup
2\cup J_2}}}b_{J_2-2}(p)b_{J_1}(p)A_I(p,p^{-s})\\& =
\sum_{\underset{J_2\subseteq\{3,4,5\}}{J_1\subseteq\{1^*\}}}b_{J_2-2}(p)b_{J_1}(p)(A_{J_1\cup
2^* \cup J_2}(p,p^{-s}) - A_{J_1\cup 2\cup J_2}(p,p^{-s}) )
\end{align*}

And finally we do the same for $\n_3(p)$ and get

\begin{align*}
W_3(p,p^{-s})&=\sum_{\underset{J_2\subseteq\{4,5\}}{\underset{J_1\subseteq\{1^*,2^*\}}
{I=J_1\cup
3^*\cup J_2}}}b_{J_2-3}(p)b_{J_1}(p)A_I(p,p^{-s})-
\sum_{\underset{J_2\subseteq\{4,5\}}{\underset{J_1\subseteq\{1^*,2^*\}}{I=J_1\cup
3\cup J_2}}}b_{J_2-3}(p)b_{J_1}(p)A_I(p,p^{-s})\\& =
\sum_{\underset{J_2\subseteq\{4,5\}}{J_1\subseteq\{1^*,2^*\}}}b_{J_2-3}(p)b_{J_1}(p)(A_{J_1\cup
3^* \cup J_2}(p,p^{-s}) - A_{J_1\cup 3\cup J_2}(p,p^{-s}) )
\end{align*}

\subsection{Igusa factors}

We observe the factor $I_n(\mathbf{U})$ appearing frequently in
the formulae. In this section we see that it is a natural part of
the zeta functions.

\begin{Def} In an expression of the form
$$W_i(X,Y)=I_{d'-i-1}(X_{d'-1},\dots,X_{i+1})E_i(X_i,Y_i)I_{i-1}(Y_{i-1},\dots,Y_1)$$
we shall call the factor $I_{d'-i-1}(\mathbf{X})$ the \emph{upper
Igusa factor} and the factor $I_{i-1}(\mathbf{Y})$ the \emph{lower
Igusa factor}. \end{Def}

Some of Voll's work has concerned the Igusa factor, for instance his formula for the normal zeta
function of the so-called Grenham groups in \cite{FE} is completely of this form. There he also
observed the existence of the upper Igusa factor in the formula of $W_1(p,p^{-s})$. The existence
of the lower Igusa factor is recorded here publicly for the first time.
However, it did appear in the
calculation of the
Segre example \cite{Segre} where the Pfaffian hypersurface was the Segre surface.

Let us now start to calculate and determine where these factors
come from.

For completeness we shall calculate  $W_0(p,p^{-s})$ here again,
and in this context using the $\mu$-function.

\begin{lemma}[Voll]
\begin{equation}
W_0(p,p^{-s})=\sum_{I\subseteq\{1,\dots,5\}}
b_I(p)A_I(p,p^{-s})=\sum_{I\subseteq\{1,\dots,5\}}b_{I}(p^{-1})
\prod_{i\in I}\frac{X_i}{1-X_i},\end{equation} where
$X_i=p^{i(10-i)-(4+i)s}$ for $i=1,\dots,5.$
\end{lemma}

\begin{proof}
Using the lemma \ref{no of lattices} we can write
\begin{align*}&\sum_{I\subseteq\{1,\dots,5\}}
b_I(p)A_I(p,p^{-s})\\&=\sum_{I\subseteq\{1,\dots,5\}}
b_I(p)\prod_{i\in I}\sum_{r_i=1}^\infty p^{4ir_i-(4+i)r_is}
p^{-\dim\flag_I}p^{\sum_{i\in I} (6-i)ir_i}\\&=
\sum_{I\subseteq\{1,\dots,5\}}b_I(p^{-1})\prod_{i\in I}\sum_{r_i=1}^\infty
p^{i(10-i)r_i-(4+i)r_is}\\&=
\sum_{I\subseteq\{1,\dots,5\}}b_I(p^{-1})\prod_{i\in I} \frac{p^{i(10-i)-(4+i)s}}{1-
p^{i(10-i)-(4+i)s}}.
\end{align*}
\end{proof}

\subsection{Extracting the upper Igusa factor}

For simplicity of notation we shall ignore the variables  $p$ and
$p^{-s}$ from the generating functions and use only $A_I$ to
denote $A_I(p,p^{-s}).$ We shall also from now on forget the
$p-$adic valuations $v_p(a_{ij})$ in the $\min$-expressions and
will write, for example, $\min\{r_1,a_{15}\}$ instead of
$\min\{r_1,v_p(a_{15})\}.$

We can
first simplify our formulae by summing out most of the variables.
We do this in the next two lemmas.

\begin{lemma}
\begin{equation}(A_{J_1 \cup i^* \cup J_2}-A_{J_1\cup i\cup J_2})=
p^{-\dim \flag_{J_2-i}} \tilde{A}_{J_2}\cdot(A_{J_1\cup i^*}-A_{J_1\cup i}),
\end{equation} where $\tilde{A}_{J_2}=\prod_{j_2\in J_2} \frac{X_{j_2}}{1-X_{j_2}}$ and $X_{j_2}=p^{(10-j_2)j_2-s(4+j_2)}$ as before.
\end{lemma}

\begin{proof}
\begin{align*}
&(A_{J_1 \cup i^* \cup J_2}-A_{J_1\cup i\cup J_2})=
\prod_{j_2\in J_2}\sum_{r_{j_2}=1}^\infty p^{4j_2r_{j_2}-s(4+j_2)r_{j_2}}
\sum_{r_i=1}^\infty p^{4ir_i-s(4+i)r_i}\cdot\\& \cdot
\prod_{j_1\in J_1}\sum_{r_{j_1}=1}^\infty p^{4j_1r_{j_1}-s(4+j_1)r_{j_1}}
\prod_{\underset{j\leq i}{i,j\in J_1\cup i\cup J_2}}
\sum_{a_{ij}=1}^{\sum_{k=i}^j r_k }\mu (\sum_{k=i}^j r_k, a_{ij})\cdot\\ &\cdot
(p^{s\min\{\sum_{j_1\in J_1} r_{j_1}+r_i, a_{ij}+r_1+\dots+r_{i-1}\}}-p^{s\min\{\sum_{j_1\in J_1} r_{j_1}, a_{ij}+r_1+\dots+r_{i-2}\}}).
\end{align*}
We can now sum up all those $\mu$-functions that depend only on
$r_{j_2} $ and from those we get
$p^{-\dim \flag_{J_2-i}}p^{\sum_{j_2\in J_2}(6-j_2-i)j_2r_{j_2}}$. Moreover we
can use the property of $\mu$ that $$\sum_{a_{ij}=1}^{r_i+\dots+r_j}
\mu(r_i+\dots+r_j,a_{ij})=p^{r_{k}+\dots +r_j}
\sum_{a_{ij}=1}^{r_i+\dots+r_{k-1}}\mu(r_i+\dots+r_{k-1},a_{ij})$$ to extract the $r_{j_2}$ from the $\mu$-part as well, as these don't appear in the $\min$-expressions. Then we get $p^{\sum_{j_2\in J_2}  ij_2r_{j_2}}$ out of the  $\mu$'s.

In conclusion, we have extracted $$p^{-\dim
\flag_{J_2-i}}p^{\sum_{j_2\in J_2}(6-j_2)j_2r_{j_2}},$$ and we can
rewrite

\begin{align*}(A_{J_1 \cup i^* \cup J_2}-A_{J_1\cup i\cup J_2})&=
\prod_{j_2\in J_2}\sum_{r_{j_2}=1}^\infty p^{4j_2r_{j_2}-s(4+j_2)r_{j_2}}
p^{-\dim \flag_{J_2-i}}p^{\sum_{j_2\in J_2}(6-j_2)j_2r_{j_2}}
(A_{J_1\cup i^*}-A_{J_1\cup i})\\&=
p^{-\dim \flag_{J_2-i}} \prod_{j_2\in J_2}\sum_{r_{j_2}=1}^\infty p^{(10-j_2)r_{j_2}-s(4+j_2)r_{j_2}} (A_{J_1\cup i^*}-A_{J_1\cup i})\\&=
p^{-\dim \flag_{J_2-i}}\tilde{A}_{J_2}\cdot(A_{J_1\cup i^*}-A_{J_1\cup i})
\end{align*}

\end{proof}

\begin{lemma}
\begin{align*}
W_i&=\sum_{\underset{J_2\subseteq\{i+1,\dots,5\}}{J_1
\subseteq\{1^*,\dots,{i-1}^*\}}}b_{J_2-i}(p)b_{J_1}(p)(A_{J_1\cup
i^* \cup J_2} - A_{J_1\cup i\cup J_2} )\\ &=I_{5-i}(X_{i+1},\dots,X_5)\cdot
 \sum_{J_1\subseteq\{1^*,\dots,{i-1}^*\}}b_{J_1}(p)
(A_{J_1\cup i^*}-A_{J_1\cup i})
\end{align*} with the poles $X_k=p^{k(10-k)-(4+k)s}$ for $k=i+1,\dots,5.$
\end{lemma}

\begin{proof}
\begin{align*}
W_i&=\sum_{\underset{J_2\subseteq\{i+1,\dots,5\}}{J_1
\subseteq\{1^*,\dots,{i-1}^*\}}}b_{J_2-i}(p)b_{J_1}(p)(A_{J_1\cup
i^* \cup J_2} - A_{J_1\cup i\cup J_2}) \\ &=
\sum_{\underset{J_2\subseteq\{i+1,\dots,5\}}{J_1
\subseteq\{1^*,\dots,{i-1}^*\}}}b_{J_2-i}(p)b_{J_1}(p)
p^{-\dim \flag_{J_2-i}}A_{J_2}\cdot(A_{J_1\cup i^*}-A_{J_1\cup i})\\&=
\sum_{J_2\subseteq\{i+1,\dots,5\}}\sum_{J_1
\subseteq\{1^*,\dots,{i-1}^*\}}b_{J_2-i}(p^{-1})b_{J_1}(p)
A_{J_2}\cdot(A_{J_1\cup i^*}-A_{J_1\cup i})\\ &=
\sum_{J_2\subseteq\{i+1,\dots,5\}}b_{J_2-i}(p^{-1})A_{J_2}
\sum_{J_1
\subseteq\{1^*,\dots,{i-1}^*\}}b_{J_1}(p)(A_{J_1\cup i^*}-A_{J_1\cup i})\\&=
I_{5-i}(X_{i+1},\dots,X_5)\cdot
\sum_{J_1
\subseteq\{1^*,\dots,{i-1}^*\}}b_{J_1}(p)(A_{J_1\cup i^*}-A_{J_1\cup i}).
\end{align*}
\end{proof}

\begin{cor}\label{essentials}
Using the above lemma we can write the rational functions
\begin{align*}
W_0(p,p^{-s})&=\sum_{I\subseteq\{1,\dots,5\}} b_I(p^{-1})
\prod_{i\in I}\frac{X_i}{1-X_i}=I_5(X_1,\dots,X_5).\\
W_1(p,p^{-s})&=I_4(X_2,\dots,X_5)\cdot(A_{1^*}-A_{1})\\
W_2(p,p^{-s})&=I_3(X_3,X_4,X_5)\cdot \sum_{J_1\subseteq\{1^*\}}b_{J_1}(p)
(A_{J_1\cup 2^*}-A_{J_1\cup 2})\\
W_3(p,p^{-s})&=I_2(X_4,X_5)\cdot
\sum_{J_1\subseteq\{1^*,2^*\}}b_{J_1}(p)(A_{J_1\cup 3^*}-A_{J_1\cup 3}),\end{align*}
where $X_i=p^{(10-i)i-(4+i)s}$ for $i=1,\dots,5.$
\end{cor}

Thus it is enough to consider the sums
$$\sum_{J_1\subseteq\{1^*,\dots,{i-1}^*\}}b_{J_1}(p)(A_{J_1\cup i^*}-A_{J_1\cup i}),$$
and that is what we shall do in the next section.

\subsection{The lower Igusa factor}

From the corollary \ref{essentials} we recall the form of $W_2(p,p^{-s})$ as
\begin{equation}
W_2(p,p^{-s})=I_3(X_3,X_4,X_5)\cdot \sum_{J_1\subseteq\{1^*\}}b_{J_1}(p)
(A_{J_1\cup 2^*}-A_{J_1\cup 2})\end{equation}

Here we need to calculate two terms, namely $(A_{2^*}-A_2)$ and
$(A_{1^*,2^*}-A_{1^*,2})$.

We shall start with the latter. This example demonstrates very
well the general idea and isn't too hard to do explicitly.

\begin{lemma}\label{idea of summations}
\begin{equation}(A_{1^*,2^*}-A_{1^*,2})=p^{-\dim\flag_{\{1\}}}\frac{Y_1}{1-Y_1}(A_{2^*}-A_2),\end{equation} where $Y_1=p^{8-3s}$.
\end{lemma}

\begin{cor}\label{W2}
\begin{equation*}\sum_{J_1\subseteq\{1^*\}}b_{J_1}(p)
(A_{J_1\cup 2^*}-A_{J_1\cup 2})=I_1(Y_1)(A_{2^*}-A_2).
\end{equation*}
\end{cor}

\begin{proof}(Proof of corollary \ref{W2})
\begin{align*}\sum_{J_1\subseteq\{1^*\}}b_{J_1}(p)
(A_{J_1\cup 2^*}-A_{J_1\cup 2})&=(1+b_{1}(p)p^{-1}\frac{Y_1}{1-Y_1})(A_{2^*}-A_2)\\&=I_1(Y_1)(A_{2^*}-A_2)\end{align*} by definition of the Igusa factor.
\end{proof}

So we have proved that $W_2(p,p^{-s})$ is of the form claimed. It
is left to prove the lemma. This is the first example of a process
I call extracting the lower Igusa factor.

\begin{proof} (Proof of lemma \ref{idea of summations})
For simplicity, let us write $T:=p^{-s}$.

%Then
%\begin{align*}A_{1^*,2^*}(p,T)&=\sum_{r_{2}=1}^\infty
%p^{8r_2}T^{6r_{2}}\sum_{r_1=1}^\infty p^{4r_3}T^{5r_3}
%\sum_{a_{15}=1}^{r_1+r_2}\mu(r_1+r_2,a_{15})\sum_{a_{24},a_{25}=1}^{r_2}
%\mu(r_2;a_{24},a_{25})
%\cdot\\ & \cdot T^{-\min\{r_1,a_{15}
%\}-\min\{r_1+r_2,a_{15},a_{25}+r_1,a_{24}+r_1\} }.
%\end{align*}
%and
%\begin{align*}A_{1^*,2}(p,T)&=\sum_{r_{2}=1}^\infty
%p^{8r_2}T^{6r_{2}}\sum_{r_1=1}^\infty p^{4r_3}T^{5r_3}
%\sum_{a_{15}=1}^{r_1+r_2}\mu(r_1+r_2,a_{15})\sum_{a_{24},a_{25}=1}^{r_2}
%\mu(r_2;a_{24},a_{25})
%\cdot\\ & \cdot T^{-2\min\{r_1,a_{15}
%\}}.
%\end{align*} These don't differ much and we can write

\begin{align*}
A_{1^*,2^*}-A_{1^*,2}= &\sum_{r_{2}=1}^\infty
p^{8r_2}T^{6r_{2}}\sum_{r_1=1}^\infty p^{4r_1}T^{5r_1}
\sum_{a_{15},a_{14},a_{13},a_{12}=1}^{r_1+r_2}\mu(r_1+r_2;a_{15},a_{14},a_{13},a_{12})\\&
\sum_{a_{11}=1}^{r_1}\mu(r_1,a_{11})
\sum_{a_{25},a_{24},a_{23},a_{22}=1}^{r_2}
\mu(r_2;a_{25},a_{24},a_{23},a_{22})
\cdot\\ & \cdot (T^{-\min\{r_1,a_{15}
\}-\min\{r_1+r_2,a_{15},a_{25}+r_1,a_{24}+r_1\} }-T^{-2\min\{r_1,a_{15}
\}} ).
\end{align*}

First we can see that some variables are independent and can be
summed separately. We shall write
$\sum_{a_{11}=1}^{r_1}\mu(r_1,a_{11})=p^{r_1-1}$ and more
importantly $$\sum_{a_{15},a_{14},a_{13},a_{12}=1}^{r_1+r_2}
\mu(r_1+r_2; a_{15},a_{14},a_{13},a_{12})=p^{3r_1}
\sum_{a_{14},a_{13},a_{12}=1}^{r_2}\mu(r_2;a_{14},a_{13},a_{12})
\sum_{a_{15}=1}^{r_1+r_2} \mu(r_1+r_2,a_{15}).$$

Denoting
$\mathbf{a}=(a_{14},a_{13},a_{12},a_{25},a_{24},a_{23},a_{22})$
and inserting the above two simplifications we can write our
summations as
\begin{align*}
A_{1^*,2^*}-A_{1^*,2}= &p^{-1}\sum_{r_{2}=1}^\infty
p^{4r_2}T^{6r_{2}}\sum_{\mathbf{a}=1}^{r_2}\mu(r_2;\mathbf{a})
\sum_{r_1=1}^\infty p^{8r_1}T^{5r_1}
\sum_{a_{15}=1}^{r_1+r_2}\mu(r_1+r_2,a_{15})\cdot\\
& \cdot (T^{-\min\{r_1,a_{15}
\}-\min\{r_1+r_2,a_{15},a_{25}+r_1,a_{24}+r_1\} }-T^{-2\min\{r_1,a_{15}
\}} ).\end{align*}

Write
\begin{align*}
B(r_2)&=\sum_{r_1=1}^\infty p^{8r_1}T^{5r_1}
\sum_{a_{15}=1}^{r_1+r_2}\mu(r_1+r_2,a_{15})\cdot\\
& \cdot (T^{-\min\{r_1,a_{15}
\}-\min\{r_1+r_2,a_{15},a_{25}+r_1,a_{24}+r_1\} }-T^{-2\min\{r_1,a_{15}
\}} )
\end{align*} and start to decompose this according $r_1=1$ and $r_1\geq 2$,
in the latter case change the variable $r'_1=1+r_1$:

\begin{align*}
B(r_2)&=p^{8}T^{5} \sum_{a_{15}=1}^{1+r_2}\mu(1+r_2,a_{15})\cdot\\
& \cdot (T^{{-1}-\min\{1+r_2,a_{15},a_{25}+1,a_{24}+1\}
}-T^{-2\min\{1,a_{15} \}}) \\& +p^8T^5 \sum_{r'_1=1}^\infty
p^{8r_1}T^{5r_1}
\sum_{a_{15}=1}^{1+r'_1+r_2}\mu(1+r'_1+r_2,a_{15})\cdot\\ & \cdot
(T^{-\min\{1+r'_1,a_{15}
\}-\min\{1+r'_1+r_2,a_{15},a_{25}+r'_1+1,a_{24}+r'_1+1\}
}-T^{-2\min\{1+r'_1,a_{15} \}} ).
\end{align*}

Now if $a_{15}=1$ all the $\min$-expressions will take value $1$
and so the subtraction will cancel all these parts. Thus we can
assume that $a_{15}\geq 2$. Moreover, we can change the range of
summation $$\sum_{a_{15}=2}^{1+r'_1}\mu(1+r'_1,a_{15})
=\sum_{a_{15}=1}^{r'_1}\mu(r'_1,a_{15})$$ and note that
$\min\{1+r'_1,1+a_{15}\}=1+\min\{r'_1,a_{15}\}$. Now we can write
the summations in the form
\begin{align*}
B(r_2)&=p^{8}T^{3} \sum_{a_{15}=1}^{r_2}\mu(r_2,a_{15})\cdot\\ &
\cdot (T^{-\min\{r_2,a_{15},a_{25},a_{24}\} }-1) \\& +p^8T^3
\sum_{r'_1=1}^\infty p^{8r_1}T^{5r_1}
\sum_{a_{15}=1}^{r'_1+r_2}\mu(r'_1+r_2,a_{15})\cdot\\ & \cdot
(T^{-\min\{r'_1,a_{15}
\}-\min\{r'_1+r_2,a_{15},a_{25}+r'_1,a_{24}+r'_1\}
}-T^{-2\min\{r'_1,a_{15} \}} ).
\end{align*}

Inserting this back to $A_{1^*,2^*}-A_{1^*,2}$ we obtain that
\begin{equation}
(1-p^8T^3)(A_{1^*,2^*}-A_{1^*,2})=p^{-1}p^8T^3(A_{2^*}-A_{2}),
\end{equation}
and thus
\begin{equation*}(A_{1^*,2^*}-A_{1^*,2})=p^{-1}
\frac{Y_1}{1-Y_1}(A_{2^*}-A_{2})\end{equation*}
where $Y_1=p^8T^3$ as claimed.

\end{proof}

Calculating $W_3(p,p^{-s})$ is very similar to $W_2(p,p^{-s})$. We
just have some more summands to consider.

\begin{prop}\label{W3}
\begin{equation}
\sum_{J_1\subseteq\{1^*,2^*\}}b_{J_1}(p)(A_{J_1\cup
3^*} - A_{J_1\cup 3})= I_2(Y_1,Y_2)\cdot (A_{3^*}-A_{3}).
\end{equation}
\end{prop}

Let us state the lemmas needed

\begin{lemma}
$$(A_{1^*\cup 3^*} - A_{1^*\cup 3})=p^{-\dim\flag_{1}}\frac{Y_1}{1-Y_1}
(A_{3^*}-A_3)$$
\end{lemma}

\begin{lemma}
$$(A_{2^*\cup 3^*} - A_{2^*\cup 3})=p^{-\dim\flag_{2}}\frac{Y_2}{1-Y_2}
(A_{3^*}-A_3)$$
\end{lemma}

\begin{lemma}\label{mixed}
$$(A_{1^*,2^*\cup
3^*} - A_{1^*,2^*\cup 3})=p^{-\dim\flag_{1,2}}\frac{Y_2}{1-Y_2}\frac{Y_1}{1-Y_1}(A_{3^*}-A_3)$$
\end{lemma}

\begin{proof}
The above three lemmas are proved in exactly the same way as lemma
\ref{idea of summations} -- in the proof of lemma \ref{mixed} we
just need to do the procedure twice.
\end{proof}

\begin{proof}(of proposition \ref{W3})
By the above lemmas
\begin{align*}\sum_{J_1\subseteq\{1^*,2^*\}}b_{J_1}(p)(A_{J_1\cup
3^*} - A_{J_1\cup 3})&= \sum_{J_1\subseteq\{1^*,2^*\}}b_{J_1}(p^{-1})\prod_{j_1\in J_1}\frac{Y_{j_1}}{1-Y_{j_1}}\\&=
I_2(Y_1,Y_2)\cdot (A_{3^*}-A_{3}).
\end{align*}
\end{proof}

\subsection{The exceptional factor and final formulae}

Finally we need to calculate the exceptional factor $E_i=(A_{i^*}-A_i)$.

\begin{prop}\label{crucial corollary}
Let $d=h(G^{ab})$ be the torsion-free rank of the abelianisation
and $d'=h(Z(G))$ the torsion-free rank of the centre. Let
$n_i=c_i+d_i$ where $n_i$ is the dimension of the space of
$i-1$-dimensional linear subspace in $\P^{d'-1}$, $d_i$ the
dimension of the Fano variety of $i-1$-dimensional linear subspace
of the Pfaffian, and $c_i$ the codimension of the same object.
Then
\begin{align*}
A_{i^*}-A_i&=
\sum_{r_i=1}^{\infty}\sum_{b_1=1}^{r_i}\sum_{b_2=1}^{r_i}\dots
\sum_{b_{n_i}=1}^{r_i}\mu(r_i;b_1,b_2,\dots,b_{n_i})p^{idr_i}T^{(d+i)r_i}
(T^{-t\min\{r_i,b_1,b_2,\dots,b_{c_i}\}}-1)\\
&=\frac{p^{di}T^{d+i-t}(1-T^t)}{(1-p^{id+d_i}
T^{d+i-t})(1-p^{di+n_i}T^{d+i})}\end{align*} where
$t=2$ in the case of points
and $t=1$ in the case of lines and higher dimensional Fano
varieties.
\end{prop}

To prove the proposition, we use the following
{\bf{crucial lemma}}.

\begin{lemma}\label{crucial lemma}
\begin{align*}&\sum_{r_i=1}^{\infty}\sum_{b_1=1}^{r_i}\sum_{b_2=1}^{r_i}\dots
\sum_{b_{c_i}=1}^{r_i}
\mu(r_i;b_1,b_2,\dots,b_{c_i})p^{(id+d_i)r_i}T^{(d+i)r_i}T^{-t\min\{a,b_1,b_2,\dots,b_{c_i}\}}\\&=\frac{p^{id+d_i}T^{d+i-t}(1-p^{id+d_i}T^{d+i})}
{(1-p^{id+d_i}T^{d+i-t})(1-p^{id+d_i+c_i}T^{d+i})}
\end{align*}
\end{lemma}

Now we can prove the proposition.

\begin{proof}(of proposition \ref{crucial corollary}.)
\begin{equation*}\begin{split}A_{i^*}-A_i&=
\sum_{r_i=1}^{\infty}\sum_{b_1=1}^{r_i}\sum_{b_2=1}^{r_i}\dots
\sum_{b_{n_i}=1}^{r_i}
\mu(r_i;b_1,b_2,\dots,b_n)p^{idr_i}T^{(d+i)r_i}T^{-t\min\{r_i,b_1,b_2,\dots,b_{c_i}\}}\\&-
\sum_{r_i=1}^{\infty}\sum_{b_1=1}^{r_i}\sum_{b_2=1}^{r_i}\dots
\sum_{b_{n_i}=1}^{r_i}
\mu(r_i;b_1,b_2,\dots,b_n)p^{idr_i}T^{(d+i)r_i}\\
\end{split}\end{equation*}

First we sum all the $b_i$'s that don't appear in the
$\min$-expression using the properties of the $\mu$-function:

\begin{equation*}\begin{split}
A_{i^*}-A_i&=
p^{-d_i}\sum_{r_i=1}^{\infty}\sum_{b_1=1}^{r_i}\sum_{b_2=1}^{r_i}\dots
\sum_{b_{c_i}=1}^{r_i}
\mu(r_i,b_1,b_2,\dots,b_{c_i})p^{(id+d_i)r_i}T^{(d+i)r_i}
T^{-t\min\{r_i;b_1,b_2,\dots,b_{c_i}\}}\\& -p^{-(d_i+c_i)}
\sum_{r_i=1}^{\infty}p^{(id+d_i+c_i)r_i}T^{(d+i)r_i}\\
\end{split}\end{equation*}

Now using lemma \ref{crucial lemma} and summation of geometric
progressions we get

\begin{align*}A_{i^*}-A_i&=
p^{-d_i}\frac{p^{id+d_i}T^{d+i-t}(1-p^{id+d_i}T^{d+i})}
{(1-p^{id+d_i}T^{d+i-t})(1-p^{id+d_i+c_i}T^{d+i})}\\&
-p^{-(d_i+c_i)}\frac{p^{id+d_i+c_i
}}{1-p^{id+c_i+d_i}T^{d+i}},\end{align*} and by a routine
calculation we get the formula given in the corollary.

\end{proof}

And finally, the proof of  lemma \ref{crucial lemma} depends on
the following combinatorial observations:

\begin{lemma}\label{shifting lemma}
\begin{align*}
&\sum_{b_1=1}^{r_i+1}\sum_{b_2=1}^{r_i+1}\dots \sum_{b_{c_i}=1}^{r_i+1}
\mu(r_i+1;b_1,b_2,\dots,b_{c_i})T^{-t\min\{r_i+1,b_1,b_2,\dots,b_{c_i}\}}\\ &=
T^{-t}\sum_{b_1=1}^{r_i}\sum_{b_2=1}^{r_i}\dots \sum_{b_{c_i}=1}^{r_i}
\mu(r_i;b_1,b_2,\dots,b_{c_i})T^{-t\min\{r_i,b_1,b_2,\dots,b_{c_i}\}}+T^{-t}p^{c_ir_i}(1-p^{-{c_i}})
\end{align*}
\end{lemma}

\begin{lemma}\label{binomial lemma}
$$p^{nk}(1-p^{-n})=\binom{n}{1}p^k(1-p^{-1})(p^{k-1})^{n-1}
+\binom{n}{2}(p^k(1-p^{-1}))^2(p^{k-1})^{n-2}
+\dots+(p^k(1-p^{-1}))^n$$
\end{lemma}

\begin{proof} Binomial theorem. \end{proof}

\begin{lemma}\label{translation lemma}
\begin{align*}
&\sum_{b_1=2}^{r_i+1}\sum_{b_2=2}^{r_i+1}\dots \sum_{b_{c_i}=2}^{r_i+1}
\mu(r_i+1;b_1,b_2,\dots,b_{c_i})T^{-t\min\{r_i+1,b_1,b_2,\dots,b_{c_i}\}}\\& =
T^{-t}\sum_{b_1=1}^{r_i}\sum_{b_2=1}^{r_i}\dots \sum_{b_{c_i}=1}^{r_i}
\mu(r_i;b_1,b_2,\dots,b_{c_i})T^{-t\min\{r_i,b_1,b_2,\dots,b_{c_i}\}}
\end{align*}
\end{lemma}

\begin{proof}
Change the range of summation and observe that
$\mu$ is invariant under this action.
\end{proof}Next we prove lemma \ref{shifting lemma}.

\begin{proof}
We split the range of summation into two parts
according to whether the min-expression takes value equal or
greater than 1. Then
\begin{align*}&\sum_{b_1=1}^{r_i+1}\dots\sum_{b_{c_i}=1}^{r_i+1}
\mu(r_i+1;b_1,\dots,b_{c_i})T^{-t\min\{r_i+1,b_1,\dots,b_{c_i}\}}
=\\&=T^{-t}(\text{number
of terms for which min is = 1})+\\& +
\sum_{b_1=2}^{r_i+1}\dots\sum_{b_{c_i}=2}^{r_i+1}\mu(r_i+1;b_1,\dots,b_{c_i})T^{-t\min\{r_i+1,b_1\dots,b_{c_i}\}}
\\&= T^{-t}\cdot(\text{no of
\{min=1\}})+T^{-t}\sum_{b_1=1}^{r_i}\dots\sum_{b_{c_i}=1}^{r_i}\mu(r_i;,b_1,\dots,b_{c_i})
T^{-t\min\{r_i,b_1,\dots,b_{c_i}\}}
\end{align*}
by lemma \ref{translation lemma}. Now all that is left is to count
the terms for which $\min=1$. As all the $b_i$ are integers
between 1 and $r_i+1$, the only way that
$\min\{r_i+1,b_1,\dots,b_{c_i}\}=1$ is when one or more of the
$b_i=1$. Note that $\mu(r_i+1,1)=p^{r_i}(1-p^{-1})$, and
$\sum_{b_i=2}^{r_i+1}\mu(r_i+1,b_i)=p^{r_i-1}$, so with these we
have multiplicity
$\binom{{c_i}}{r}(p^{r}(1-p^{-1}))^m(p^{{r_i}-1})^{{c_i}-m}$ when
exactly $m$ of the $b_i$ are equal to one. Now putting these
together and using lemma \ref{binomial lemma} we have the
multiplicity of $T^{-t}$ $p^{c_ir_i}(1-p^{-c_i})$ as claimed. This
proves lemma \ref{shifting lemma}.
\end{proof}

And then finally lemma \ref{crucial lemma}:

\begin{proof}(of lemma \ref{crucial lemma})
Write $$A'_{i^*}=\sum_{r_i=1}^{\infty}
p^{(id+d_i)r_i}T^{(d+i)r_i}\sum_{b_1=1}^{r_i}\sum_{b_2=1}^{r_i}\dots
\sum_{b_{c_i}=1}^{r_i}
\mu(r_i;b_1,b_2,\dots,b_{c_i})T^{-t\min\{r_i,b_1,b_2,\dots,b_{c_i}\}}$$
We split $A'_{i^*}$ into pieces when $r_i=1$ and $r_i\geq 2$ and
in the latter we change the variable $r_i=k+1$.

\begin{align*}
A'_{i^*}&=p^{id+d_i}T^{d+i-t}+\\&+p^{id+d_i}T^{d+i}\sum_{k=1}^{\infty}p^{(id+d_i)k}T^{(d+i)k}\sum_{b_1=1}^{k+1}\dots
\sum_{b_{c_i}=1}^{k+1}
\mu(k+1;b_1,b_2,\dots,b_{c_i})T^{-t\min\{k+1,b_1,b_2,\dots,b_{c_i}\}}\\
&=p^{id+d_i}T^{d+i-t}+\\&+ p^{id+d_i}T^{d+i-t}\sum_{k=1}^{\infty}
p^{(id+d_i)k}T^{(d+i)k} \sum_{b_1=1}^k\dots \sum_{b_{c_i}=1}^k
\mu(k,b_1,b_2,\dots,b_{c_i})T^{-t\min\{k;b_1,b_2,\dots,b_{c_i}\}}\\
&+(1-p^{- {c_i}})p^{id+d_i}T^{d+i-t}\sum_{k=1}^\infty
p^{(id+d_i+c_i)k}T^{(d+i)k},
\end{align*}
using lemma \ref{shifting lemma}. Now
\begin{align*} (1-p^{(id+d_i)}T^{d+i-t})A'_{i^*}&=p^{id+d_i}T^{d+i-t}
+(1-p^{-{c_i}})p^{id+d_i}T^{d+i-t}\sum_{k=1}^\infty
p^{(id+d_i+c_i)k}T^{(d+i)k}
\\&=p^{id+d_i}T^{d+i-t}+(1-p^{-c_i})p^{id+d_i}T^{d+i-t}
\frac{p^{id+d_i+c_i}T^{d+i}}{1-p^{id+d_i+c_i}T^{d+i}},
\end{align*}
just summing up the geometric progressions. Now rearranging our
rational functions we get
$$A'_{i^*}=\frac{p^{id+d_i}T^{d+i-t}(1-p^{id+d_i}T^{d+i})}{(1-p^{id+d_i}T^{d+i-t})(1-p^{id+d_i+c_i}T^{d+i})}$$
as claimed.
\end{proof}

\section{Functional equation}\label{fe}

%Recall that the Hasse-Weil zeta function has a functional equation
%if the variety is smooth and absolutely irreducible.

The functional equation in this example is now an easy corollary
of the general shape of the zeta function. It is conjectured in
\cite{Ennui} Conjecture 5.47 that all uniform zeta functions of
class two nilpotent groups satisfy the functional equation
$$\zeta_{G,p}^\ns(s)|_{p\mapsto p^{-1}} =
(-1)^{d+d'}p^{\binom{d+d'}{2}-(2d+d')s} \cdot
\zeta_{G,p}^\ns(s).$$

Recall from \cite{Amer} the observation that
$$\zeta_{G,p}^\ns(s)|_{p\mapsto p^{-1}} =
(-1)^{d+d'}p^{\binom{d+d'}{2}-(2d+d')s} \cdot
\zeta_{G,p}^\ns(s)\iff$$ $$A(p,T)|_{p\mapsto
p^{-1}}=(-1)^{d'-1}p^{\binom{d'}{2}}\cdot A(p,T).$$

We can also observe that in our decomposition of $$A(p,T)=
W_0(p,T) + \sum_{i=1}^3 \n_i(p) W_i(p, T)$$ it is enough to show
that each of the summands satisfies the same functional equation
as $A(p,T)$ satisfies. In this example $d=4$ and $d'=6$, but we do
this in general, since the functional equation doesn't depend on
the actual numerical values $d$ and $d'$ take.

It is known at least from \cite{FE} and can in fact be explicitly observed that
$$I_k(\mathbf{U})|_{U_i\mapsto
U_i^{-1}} = (-1)^{k} p^{\binom{k+1}{2}} I_k(\mathbf{U}),$$ where
$U_i=p^{a_i-b_is}.$

So in each of
$$W_i(X,Y)=I_{d'-i-1}(X_{d'-1},\dots,X_{i+1})E_i(X_i,Y_i)I_{i-1}(Y_{i-1},\dots,Y_1)$$
we have a functional equations from the Igusa factors as
$(-1)^{d'-i-1} p^{\binom{d'-i}{2}}$ and $(-1)^{i-1}
p^{\binom{i}{2}},$ and it is left to see what is the functional
equation on $E_i$. But this is now easy since we have the explicit
formula $$E_i(X_i,Y_i)=
\frac{p^{-d_i}Y_i-p^{-n_i}X_i}{(1-X_i)(1-Y_i)},$$ and we see that
$$E_i(X_i,Y_i)|_{\mapsto p^{-1}} = p^{n_i+d_i} E_i(X_i,Y_i).$$
Here $n_i=i(d'-i)$ and $d_i=$ the dimension of the Fano variety of
$(i-1)$-dimensional linear subspaces on $\mathfrak{P}$.

If we know observe that for $i=1,2,3$ the points on the Fano
varieties have functional equations $\n_i(p)|_{p\mapsto p^{-1}}=
p^{-d_i} \n_i(p)$, since these $\n_i(p)$ are polynomials in $p$
with nice symmetry properties.

So putting all of these together $$W_i(p,T)|_{p\mapsto p^{-1}}=
(-1)^{d'-i-1+(i-1)}p^{\binom{d'-i}{2}+n_i+\binom{i}{2}}
W_i(p,T),$$ and as $\binom{d'-i}{2}+n_i+\binom{i}{2}=
\binom{d'}{2}$ we have that $$W_i(p,T)|_{p\mapsto p^{-1}}=
(-1)^{d'}p^{\binom{d'}{2}} W_i(p,T)$$ as required. This proves the
functional equation.


\begin{thebibliography}{99}


\bibitem{elliptic} M.P.F. du Sautoy, A nilpotent group and its elliptic
curve: non-uniformity of local zeta functions of groups, Israel J.
of Math. {\bf 126} (2001), 269 -- 288.


\bibitem{IHES} M.P.F. du Sautoy, Counting finite $p$-groups and
nilpotent groups, Inst. Hautes Études Scientifiques, Publ. Math.,
{\bf 92} (2000), 63 -- 112.


\bibitem{crelle} M.P.F. du Sautoy, Counting subgroups in nilpotent
groups and points on elliptic curves, J. Reine. Angew. Math. {\bf
549} (2002), 1 -- 21.


\bibitem{Ennui} M.P.F. du Sautoy, Zeta functions of groups: The
quest for order versus the flight from ennui, Groups St Andrews
2001 - in Oxford, Volume 1, CUP 2003.



\bibitem{Annals} M.P.F. du Sautoy and F.J. Grunewald, Analytic
properties of zeta functions and subgroup growth, Annals of Math.
{\bf 152} (2000), 793 -- 833.

\bibitem{GSS} F.J. Grunewald, D. Segal and G.C. Smith, Subgroups
of finite index in nilpotent groups, Invent. Math. {\bf 93}
(1988), 185 -- 223.

\bibitem{Harris} J. Harris, \emph{Algebraic Geometry: A First Course}, Springer-Verlag,
New York, 1992.

\bibitem{Hirschfeld} J.W.P. Hirschfeld and J.A. Thas,
\emph{General Galois Geometries}, Clarendon Press, Oxford, 1991.


%\bibitem{growth} A. Lubotzky and D. Segal, \emph{Subgroup growth},
%Progress in Mathematics {\bf 212}, Birkh\"auser, 2003.

%\bibitem{transfer} P.M. Paajanen, Zeta functions of groups and
%their abscissae of convergence, Transfer thesis 2003, Oxford.

\bibitem{Abscissa} P.M. Paajanen, On the degree of polynomial subgroup growth in class 2
nilpotent groups, preprint.

\bibitem{smooth} P.M. Paajanen, Structure and functional equations of normal
zeta functions of class
two nilpotent groups, in preparation.

\bibitem{Segre} P.M. Paajanen, The Segre example -- lines on Pfaffian, a chapter of DPhil thesis,
Oxford, in preparation.

\bibitem{Amer} C. Voll, Zeta functions of groups and enumeration
in Bruhat-Tits buildings, Amer. J. Math. {\bf 126} (2004), no. 5, 1005 -- 1032.

\bibitem{FE} C. Voll, Functional equations for local normal zeta
functions of nilpotent groups, Geom. Funct. Anal., to appear.

%\bibitem{PhD} C. Voll, Zeta functions of groups and enumeration in
%Bruhat-Tits buildings, Ph.D. Thesis 2002 Cambridge.

\end{thebibliography}
\end{document}